\documentclass[12pt,twoside]{article}
\usepackage{amsmath,amsfonts,amssymb,a4,enumerate,epsfig,array,theorem,psfrag}

\newtheorem{theo}{Theorem}

\newtheorem{prop}{Proposition}[section]

\newtheorem{lemma}[prop]{Lemma}

\newcommand{\xx}{{\mathbf x}}
\newcommand{\ff}{{\mathbf f}}
\newcommand{\bfg}{{\mathbf g}}
\newcommand{\bfh}{{\mathbf h}}
\newcommand{\bfn}{{\mathbf n}}

\newcommand{\ZZ}{{\mathbb{Z}}}

\newcommand{\RR}{{\mathbb{R}}}

\newcommand{\Ss}{{\mathbb{S}}}

\newcommand{\DD}{{\mathbb{D}}}

\newcommand{\HH}{{\mathbb{H}}}

\newcommand{\cL}{{\cal L}}
\newcommand{\cI}{{\cal I}}
\newcommand{\cF}{{\cal F}}
\newcommand{\cA}{{\cal A}}
\newcommand{\cC}{{\cal C}}

\newcommand{\fF}{{\mathfrak F}}

\newcommand{\nobf}{\noindent\bf}

\def\qed{\unskip\nobreak\hfil\penalty50\hskip1.75em\null\nobreak\hfil
$\blacksquare$ {\parfillskip=0pt \finalhyphendemerits=0 \par}\goodbreak}

\begin{document}
\title{The homotopy and cohomology of spaces \\
of locally convex curves in the sphere --- I}
\author{Nicolau C. Saldanha}
\maketitle

\begin{abstract}
A smooth curve $\gamma: [0,1] \to \Ss^2$ is
locally convex if its geodesic curvature is positive at every point.
J.~A.~Little showed that the space
of all locally positive curves $\gamma$
with $\gamma(0) = \gamma(1) = e_1$ and $\gamma'(0) = \gamma'(1) = e_2$
has three connected components $\cL_{-1,c}$, $\cL_{+1}$, $\cL_{-1,n}$.
The space $\cL_{-1,c}$ is known to be contractible
but the topology of the other two connected components is not well understood.
We study the homotopy and cohomology of these spaces.
In particular, for $\cL_{-1} = \cL_{-1,c} \sqcup \cL_{-1,n}$,
we show that $\dim H^{2k}(\cL_{(-1)^{k}}, \RR) \ge 1$,
that $\dim H^{2k}(\cL_{(-1)^{(k+1)}}, \RR) \ge 2$,
that $\pi_2(\cL_{+1})$ contains a copy of $\ZZ^2$
and that $\pi_{2k}(\cL_{(-1)^{(k+1)}})$ contains a copy of $\ZZ$.
\end{abstract}

\section{Introduction}

\footnote{2000 {\em Mathematics Subject Classification}.
Primary 57N65, 53C42; Secondary 34B05.
{\em Keywords and phrases} Convex curves,
topology in infinite dimension, 
periodic solutions of linear ODEs.}

A curve $\gamma: [0,1] \to \Ss^2$ is called {\sl locally convex}
if its geodesic curvature is always positive,
or, equivalently, if $\det(\gamma(t), \gamma'(t), \gamma''(t)) > 0$ for all $t$.
Let $\cL_I$ be the space of all locally convex curves $\gamma$
with $\gamma(0) = \gamma(1) = e_1$ and $\gamma'(0) = \gamma'(1) = e_2$.
The topology in this space of curves will be given by
the Sobolev metric $H^2$:
this has the minor technical advantages of making $\cL_I$
a Hilbert manifold
and of allowing for jump discontinuities in $\gamma''$ in the constructions.
The choice of metric actually makes very little difference:
since it is easy to uniformly smoothen out a curve while
keeping its geodesic curvature positive
we might just as well work with the $C^\infty$ topology,
or with $C^k$ for some $k \ge 2$.

\begin{figure}[ht]
\begin{center}
\epsfig{height=30mm,file=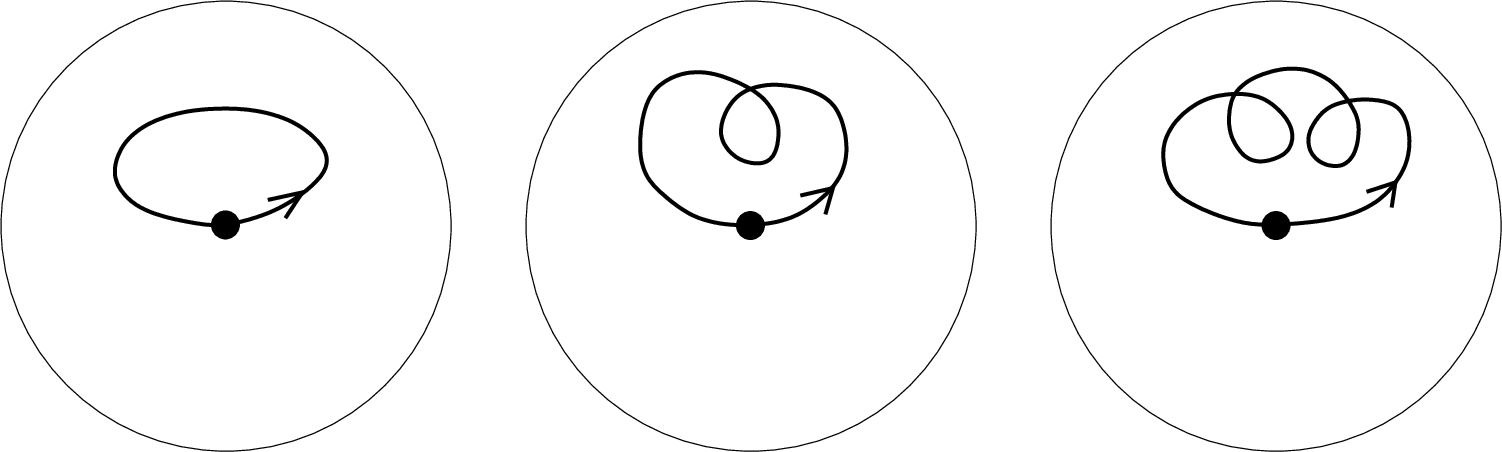}
\end{center}
\caption{Curves in $\cL_{-1,c}$, $\cL_{+1}$ and $\cL_{-1,n}$.}
\label{fig:3comp}
\end{figure}

J.~A.~Little \cite{Little} showed that $\cL_I$
has three connected components $\cL_{-1,c}$, $\cL_{+1}$, $\cL_{-1,n}$:
we call these the Little spaces.
Figure \ref{fig:3comp} shows examples of curves in
$\cL_{-1,c}$, $\cL_{+1}$ and $\cL_{-1,n}$, respectively.
The connected component $\cL_{-1,c}$ consists of the simple curves in $\cL_I$:
the space $\cL_{-1,c}$ is known to be contractible (\cite{ShapiroM}, Lemma 5).
The topology of these and related spaces has been discussed, among others,
by B. Shapiro, M. Shapiro and Khesin (\cite{Shapiro2}, \cite{SK})
but the topology of the Little spaces is still not well understood.
This series of papers is to present new results concerning
the homotopy and cohomology of the Little spaces.
A more ambitious aim would be to determine the homotopy type
of these spaces (which we hope to accomplish in \cite{Saldanha3}).

Let $\cI_I \supset \cL_I$ be the space
of immersions $\gamma: [0,1] \to \Ss^2$, $\gamma'(t) \ne 0$,
$\gamma(0) = \gamma(1) = e_1$, $\gamma'(0) = \gamma'(1) = e_2$.
For each $\gamma \in \cI_I$, consider its Frenet frame
$\fF_\gamma: [0,1] \to SO(3)$ defined by
\[ \begin{pmatrix} \gamma(t) & \gamma'(t) & \gamma''(t) \end{pmatrix} =
\fF_\gamma(t) R(t), \]
$R(t)$ being an upper triangular matrix with positive diagonal
(the left hand side is the $3\times3$ matrix with columns
$\gamma(t)$, $\gamma'(t)$ and $\gamma''(t)$).
The universal (double) cover of $SO(3)$ is $\Ss^3 \subset \HH$,
the group of quaternions of absolute value $1$;
let $\Pi: \Ss^3 \to SO(3)$ be the canonical projection.
The curve $\fF_\gamma$ can be lifted to define
$\tilde\fF_\gamma: [0, 1] \to \Ss^3$ with
$\tilde\fF_\gamma(0) = 1$, $\Pi \circ \tilde\fF_\gamma = \fF_\gamma$.
The value of $\tilde\fF_\gamma(1)$ defines the two connected components
$\cI_{\pm 1}$ of $\cI_I$:
$\gamma \in \cI_{+1}$ if and only if $\tilde\fF_\gamma(1) = 1$.
Notice that if $\gamma$ is a simple curve in $\cI_I$ then 
$\tilde\fF_\gamma(1) = -1$ and therefore $\gamma \in \cI_{-1}$.
We have $\cL_{+1} = \cI_{+1} \cap \cL_I$
and $\cL_{-1} = \cL_{-1.c} \sqcup \cL_{-1,n} = \cI_{-1} \cap \cL_I$.

Let $\Omega \Ss^3$ (resp. $\Omega_-\Ss^3$) be the set of continuous curves
$\alpha: [0,1] \to \Ss^3$, $\alpha(0) = \alpha(1) = 1$
(resp. $\alpha(0) = 1$, $\alpha(1) = -1$). 
These two spaces are easily seen to be homeomorphic 
and shall from now on be identified;
$\Omega\Ss^3$ is a well understood space:
we have $H^\ast(\Omega\Ss^3,\RR) = \RR[\xx]$
where $\xx \in H^2$ satisfies $\xx^n \ne 0$ for all positive $n$
(see, for instance, \cite{BT}).
The previous paragraph defines maps
$\tilde\fF: \cI_{\pm 1} \to \Omega\Ss^3$.
It is a well-known fact that these two maps are homotopy equivalences;
this follows from the Hirsch-Smale Theorem
(\cite{Morse}, \cite{Hirsch}, \cite{Smale}).
As we shall see,
the inclusions $i: \cL_{\pm 1} \to \cI_{\pm 1}$
are not homotopy equivalences.

\begin{theo}
\label{theo:homotosur}
For any compact space $K$ and any function
$f: K \to \cI_{\pm 1}$ there exists $g: K \to \cL_{\pm 1}$
and a homotopy $H: [0,1] \times K \to \cI_{\pm 1}$
with $H(0,k) = f(k)$ and $H(1,k) = g(k)$ for all $k \in K$.
\end{theo}

The maps $i: \cL_{\pm 1} \to \cI_{\pm 1}$
therefore induce surjective maps
$\pi_k(\cL_{\pm 1}) \to \pi_k(\cI_{\pm 1})$.
In fact, since $\cL_{\pm 1}$ and $\cI_{\pm 1}$ have the
homotopy type of CW complexes with finitely many cells
per dimension, the inclusions are homotopically surjective
but we skip the details.

In particular, $\cL_{\pm 1}$ is not homotopically equivalent
to a finite CW-complex.
Theorem 2 in \cite{SaSha} is a similar result for arbitrary dimension.
From Theorem \ref{theo:homotosur} we write
$\RR[\xx] \subseteq H^\ast(\cL_{\pm 1};\RR)$.
The main result of this paper implies that $\RR[\xx]$ has infinite codimension
as a subspace of $H^\ast(\cL_{\pm 1})$.

Let $\cL$ be the contractible space of all locally convex curves $\gamma$
with $\gamma(0) = e_1$ and $\gamma'(0) = e_2$
and define $\phi: \cL \to \Ss^3$ by
$\phi(\gamma) = \tilde\fF_\gamma(1)$
where $\tilde\fF_\gamma: [0,1] \to \Ss^3$
is defined via Frenet frames as above.
It is natural to conjecture that $\phi: \cL \to \Ss^3$
is somehow similar to a fibration.

This is not the case: we prove that the map $\phi$
does not satisfy the homotopy lifting property.
Let $X = [0,1]^2$.
There exist maps $\bfh: X \times [0,1] \to \Ss^3$ and
$\tilde\bfh_0: X \times \{0\} \to \cL$ such that
$\phi \circ \tilde\bfh_0 = \bfh|_{X \times \{0\}}$
and there exists no map $\tilde\bfh: X \times [0,1] \to \cL$
with $\phi \circ \tilde\bfh = \bfh$.
The maps $\bfh$ and $\tilde\bfh_0$ will be constructed
explicitly in Section \ref{sect:homolift}.

\begin{figure}[ht]
\begin{center}
\epsfig{height=30mm,file=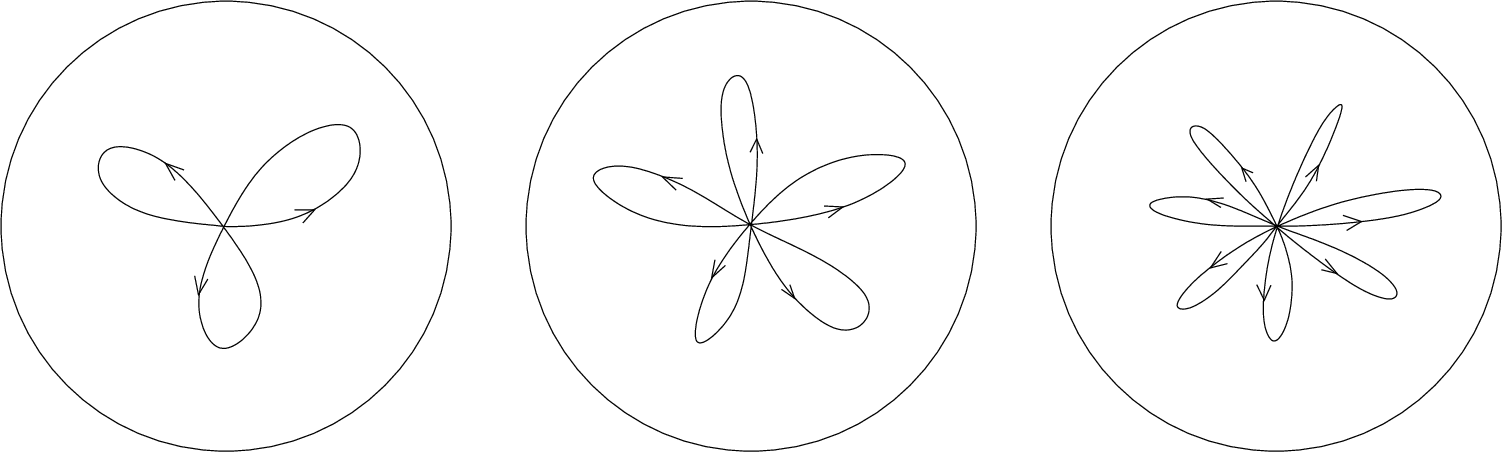}
\end{center}
\caption{Examples of flowers with $3$, $5$ and $7$ petals.}
\label{fig:flowers}
\end{figure}

A curve $\gamma \in \cL_{(-1)^{(k+1)}}$ is a {\sl flower} of $2k+1$ petals if
there exist $0 = t_0 < t_1 < t_2 < \cdots < t_{2k} < t_{2k+1} = 1$ and
$0 = \theta_0 < \theta_1 < \theta_2 < \cdots < \theta_{2k} < \theta_{2k+1}
= \theta_M$ such that:
\begin{enumerate}
\item{$\gamma(t_1) = \gamma(t_2) = \cdots = \gamma(t_{2k}) = e_1$;}
\item{the only self-intersections of the curve $\gamma$
are of the form $\gamma(t_i) = \gamma(t_j)$;}
\item{the argument of $(x_{i,2}, x_{i,3})$ is $\theta_i$,
where $(0, x_{i,2}, x_{i,3}) = (-1)^i \gamma'(t_i)$.}
\end{enumerate}
As a somewhat degenerate case,
a flower of $1$ petal is a simple locally convex curve.
Figure \ref{fig:flowers} shows examples of flowers.

For $k > 0$, let $\cF_{2k} \subset \cL_{(-1)^{(k+1)}}$ be the set of
flowers of $2k+1$ petals.
As we shall see in Lemma \ref{lemma:flower},
the set $\cF_{2k}$ is closed and a submanifold of codimension $2k$.
Furthermore, the normal bundle to $\cF_{2k}$ in $\cL_{(-1)^{(k+1)}}$ is trivial.
Thus, intersection with $\cF_{2k}$ defines an element
$\ff_{2k} \in H^{2k}(\cL_{(-1)^{(k+1)}})$
with $(\ff_{2k})^2 = 0$.

We shall construct maps $\bfg_{2k}: \Ss^{2k} \to \cL_{(-1)^{(k+1)}}$
which are homotopic to a constant in $\cI_{(-1)^{(k+1)}}$
but which satisfy $\ff_{2k}\bfg_{2k} = 1$,
thus proving that both $\ff_{2k} \in H^{2k}(\cL_{\pm 1})$
and $\bfg_{2k} \in \pi_{2k}(\cL_{\pm 1})$ are nontrivial.
This establishes our main result.

\begin{theo}
\label{theo:Gcohomology}
Let $k \ge 1$. Then
$\dim H^{2k}(\cL_{(-1)^{(k+1)}}, \RR) \ge 2$
and $\pi_{2k}(\cL_{(-1)^{(k+1)}})$ contains a copy of $\ZZ$.
\end{theo}

Notice that $\cL_{per}$, the set of all $1$-periodic locally convex
curves $\tilde\gamma: \RR \to \Ss^2$, is homeomorphic to $SO(3) \times \cL_I$:
define $\Psi: \cL_{per} \to SO(3) \times \cL_I$ by
$\Psi(\tilde\gamma) =
(\fF_{\tilde\gamma}(0), (\fF_{\tilde\gamma}(0))^{-1} \tilde\gamma|_{[0,1]})$.
We usually prefer to work in $\cL_I$ but sometimes move to $\cL_{per}$.

In Section \ref{sect:convex}, we give a brief sketch of Little's Theorem
and present the concept of convex curves.
Section \ref{sect:homotosur} is dedicated to
Theorem \ref{theo:homotosur}.
In Section \ref{sect:flower}, we prove the basic facts about
the set of flowers.
Section \ref{sect:k=1} contains the construction of the map
$\bfg_2: \Ss^2 \to \cL_{+1}$ and
the proof of Theorem \ref{theo:Gcohomology} for $k = 1$.
The construction of the maps
$\bfh: X \times [0,1] \to \Ss^3$ and
$\tilde\bfh_0: X \times \{0\} \to \cL$
are presented in Section \ref{sect:homolift}.
Finally, in Section \ref{sect:Gcohomology}, we construct the maps
$\bfg_{2k}: \Ss^{2k} \to \cL_{(-1)^{(k+1)}}$ and prove
Theorem \ref{theo:Gcohomology} for $k > 1$.
Section \ref{sect:conclusion} contains a few final remarks.

In the second paper of this series (\cite{S2}) we prove that
the connected components of $\cL_I$ are simply connected and compute
the groups $H^2(\cL_{\pm 1}; \ZZ)$.

This work was motivated by an attempt to extend to
ordinary differential equations of order $3$
some of our results with
Dan Burghelea and Carlos Tomei (\cite{BST1}, \cite{BST2}).
Consider the differential equation of order 3:
\[ u'''(t) + h_1(t) u'(t) + h_0(t) u(t) = 0, \quad t \in [0,1]; \]
the set of pairs of potentials $(h_0,h_1)$ for which the equation
admits 3 linearly independent periodic solutions is homotopically
equivalent to $\cL_I$ (\cite{ST}).
The motivation of B. Shapiro and M. Shapiro for studying these spaces
is similar.


The author would like to thank Dan Burghelea and
Boris Shapiro for helpful conversations.
The author acknowledges the hospitality of
The Mathematics Department of The Ohio State University
during the winter quarters of 2004 and 2009
and the support of CNPq, Capes and Faperj (Brazil).

\section{Convex curves and Little's Theorem}

\label{sect:convex}

In this section we give a brief review of Little's argument (\cite{Little}).

Given an interval $I$, a smooth immersion $\gamma: I \to \Ss^2$ 
and $t \in I$, let $\bfn(t)$ be the unit normal vector
$\fF_\gamma(t) e_3 = \gamma(t) \times \gamma'(t)/|\gamma'(t)|$.
Given $t_0 \in I$, the function
$\eta_{t_0}(t) = \langle \gamma(t), \bfn(t_0) \rangle$
satisfies $\eta_{t_0}(t_0) = \eta'_{t_0}(t_0) = 0$.
The curve $\gamma$ is locally convex near $t_0$
if and only if $\eta''_{t_0}(t_0) > 0$.
A locally convex curve is \textit{convex} if
$\eta_{t_0}(t) > 0$ for all $t \ne t_0$.
In other words, a convex curve is contained in one of the half spaces
defined by the plane orthogonal to $\bfn(t_0)$.

Part of Little's Theorem is that the set $\cL_{-1,c}$
of simple locally compact curves is a connected component of $\cL_I$:
Little proves that simple closed locally convex curves
are convex (see also \cite{ShapiroM}). 
We shall often use this fact.

\begin{figure}[ht]
\begin{center}
\epsfig{height=35mm,file=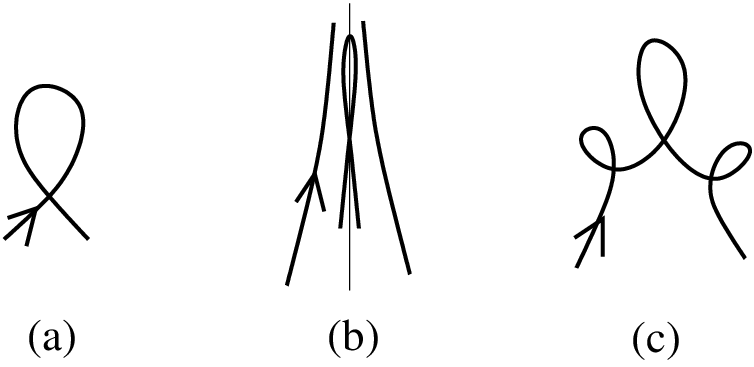}
\end{center}
\caption{Using a loop to produce two more loops.}
\label{fig:addloopl}
\end{figure}

The other part of Little's Theorem is that,
once convex curves have been removed,
the sets $\cL_{+1}$ and $\cL_{-1,n}$ are path connected.
The fundamental construction here is that
if the curve $\gamma$ has a loop,
we can add a pair of loops as in Figure \ref{fig:addloopl}:
from (a) to (b), the loop moves one full turn along a geodesic 
and from (b) to (c) the large loops are shrunk.
By repeating this procedure, we may add a large number of loops
which can then be spread along the curve.
The curve can then be deformed and, thanks to the loops,
remain locally convex.
This part will be explained in greater detail in the next section.

\section{Proof of Theorem \ref{theo:homotosur}}

\label{sect:homotosur}

\begin{figure}[ht]
\begin{center}
\epsfig{height=9mm,file=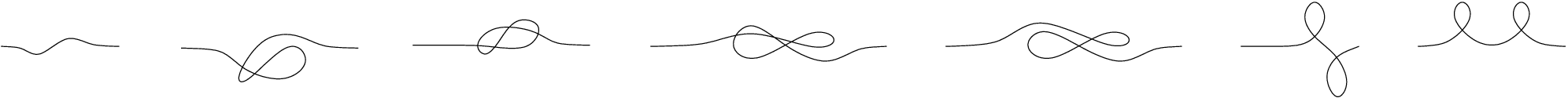}
\end{center}
\caption{How to add two small loops to a curve in $\cI_{\pm 1}$.}
\label{fig:addloopi}
\end{figure}

First notice that in $\cI_{\pm 1}$
it is easy to introduce a pair of loops at any point of the curve:
the process is illustrated in Figure \ref{fig:addloopi};
in the final step one of the loops becomes big,
goes around the sphere and shrinks again.

\begin{figure}[ht]
\begin{center}
\epsfig{height=35mm,file=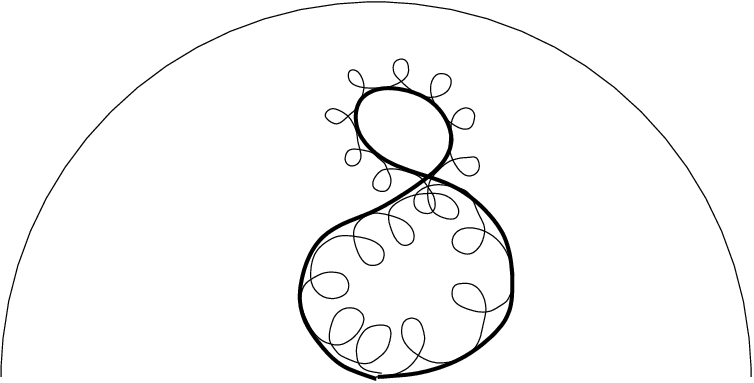}
\end{center}
\caption{Curves $\gamma_0 \in \cI_{\pm 1}$ and $\gamma_1 \in \cL_{\pm 1}$.}
\label{fig:addloop}
\end{figure}

A function $f: K \to \cI_{\pm 1}$ can be thought of as a family of curves.
We can uniformly perform the above construction several times
along all curves of the family.
Given a curve $\gamma_0$, we construct a family of curves ending
in a curve $\gamma_1$ with many positively oriented loops
as in Figure \ref{fig:addloop}.
If the number of loops is sufficiently large and the loops are tight enough,
the curve $\gamma_1$ will be locally convex.
We have therefore constructed a homotopy $H: [0,1] \times K \to \cI_{\pm 1}$
with $H(0,k) = f(k)$ and $H(1,k) \in \cL_{\pm 1}$ for all $k \in K$,
as required.

We now present a more rigorous version of this argument.
Let $\cC_0$ be the circle with diameter $e_1e_3$,
parametrized by $\nu_1 \in \cL_I$,
\[ \nu_1(t) = \left( \frac{1 + \cos(2 \pi t)}{2},
\frac{\sqrt{2}}{2} \sin(2 \pi t),
\frac{1 - \cos(2 \pi t)}{2} \right). \]
For positive $n$, let $\nu_n(t) = \nu_1(nt)$
so that $\nu_1 \in \cL_{-1,c}$ and, for $n > 1$,
$\nu_n \in \cL_{(-1)^n}$.

For $\gamma_1 \in \cI_{\sigma_1}$, $\gamma_2 \in \cI_{\sigma_2}$,
$\sigma_i \in \{+1,-1\}$,
let $\gamma_1 \ast \gamma_2 \in \cI_{\sigma_1 \sigma_2}$
be defined by
\[ (\gamma_1 \ast \gamma_2)(t) = \begin{cases}
\gamma_1(2t), & 0 \le t \le 1/2, \\
\gamma_2(2t-1), & 1/2 \le t \le 1. \end{cases} \]
Notice that if $\gamma_1, \gamma_2 \in \cL_I$
then $\gamma_1 \ast \gamma_2 \in \cL_I$.
For $f: K \to \cI_I$, let $\nu_n \ast f: K \to \cI_I$
be defined by $(\nu_n \ast f)(p) = \nu_n \ast (f(p))$.
The observation in Figure \ref{fig:addloopi} can be translated
as the following lemma, whose straightforward proof will be omitted.

\begin{lemma}
\label{lemma:addloopi}
Let $K$ be a compact set and 
let $f: K \to \cI_{\pm 1}$ a continuous function.
Then $f$ and $\nu_2 \ast f$ are homotopic.
\end{lemma}

We now need a construction corresponding to adding loops
along the curve, as in Figure \ref{fig:addloop}.
For $\gamma \in \cI_{\pm}$ and $n > 0$,
define $(F_{n}(\gamma))(t) = \fF_{\gamma}(t) \nu_{n}(t)$.

\begin{lemma}
\label{lemma:YtoX}
Let $K$ be a compact set and
let $f: K \to \cI_{\pm 1}$ a continuous function.
Then, for sufficiently large $n$,
$F_{2n} \circ f$ is a function from $K$ to $\cL_{\pm 1}$.
\end{lemma}

{\nobf Proof:}
Let $C > 1$ be a constant such that $|(\fF_\gamma)'(t)| < C$
and $|(\fF_\gamma)''(t)| < C$ for any $\gamma = f(k)$, $k \in K$
and for any $t \in [0,1]$.
Let $\epsilon > 0$ be such that if $|v_1 - \nu'_1(t)| < \epsilon$
and $|v_2 - \nu''_1(t)| < \epsilon$ then $\det(\nu_1(t), v_1, v_2) > 0$.
Take $n > 20C/\epsilon$.

For $\gamma = f(k)$, write
\[ \tilde\gamma(t) = (F_{2n}\gamma)(t) =
\fF_\gamma(t) \nu_{2n}(t) = \fF_\gamma(t) \nu_1(2nt) \]
so that
\begin{align}
\tilde\gamma'(t) &= \fF_\gamma'(t) \nu_1(2nt) +
2n \fF_\gamma(t) \nu_1'(2nt) \notag\\
\tilde\gamma''(t) &= \fF_\gamma''(t) \nu_1(2nt) +
4n \fF_\gamma'(t) \nu_1'(2nt) +
4n^2 \fF_\gamma(t) \nu_1''(2nt) \notag
\end{align}
and therefore, after a few manipulations,
\[ \left| \frac{\tilde\gamma'(t)}{2n} - \fF_\gamma(t)\nu_1'(2nt) \right|
< \epsilon, \quad
\left| \frac{\tilde\gamma''(t)}{4n^2} - \fF_\gamma(t)\nu_1''(2nt) \right|
< \epsilon \]
or, equivalently,
\[ \left| \frac{(\fF_\gamma(t))^{-1}\tilde\gamma'(t)}{2n} -
\nu_1'(2nt) \right| < \epsilon, \quad
\left| \frac{(\fF_\gamma(t))^{-1}\tilde\gamma''(t)}{4n^2} -
\nu_1''(2nt) \right| < \epsilon. \]
It follows that
\[ \det\left(\nu_\theta(2nt),\frac{(\fF_\gamma(t))^{-1}\tilde\gamma'(t)}{2n},
\frac{(\fF_\gamma(t))^{-1}\tilde\gamma''(t)}{4n^2} \right) > 0 \]
and therefore that
$\det(\tilde\gamma(t), \tilde\gamma'(t), \tilde\gamma''(t)) > 0$,
which is what we needed.
\qed

Theorem \ref{theo:homotosur} now follows directly from the next lemma.

\begin{lemma}
\label{lemma:Fandp}
Let $K$ be a compact set, $f: K \to \cI_{\pm 1}$.
Then, for sufficiently large $n$,
the image of $F_{2n} \circ f$ is contained in $\cL_{\pm 1}$ and
there exists $H: [0,1] \times K \to \cI_{\pm 1}$ such that
$H(0,\cdot) = f$ and $H(1,\cdot) = F_n \circ f$.
\end{lemma}

{\nobf Proof: }
We know from Lemma \ref{lemma:addloopi}
that $f$ if homotopic to $\nu_{2n} \ast f$.
All we have to do is construct a homotopy between $F_{2n} \circ f$
and $\nu_{2n} \ast f$.
Intuitively, this is done by pushing the loops towards $t = 0$.
More precisely, if $\gamma = f(k)$, $k \in K$, let
\[ H_1(s,k)(t) = \begin{cases} \nu_{2n}(t), &t \le s/2,\\
\fF_\gamma((2t-s)/(2-s)) \nu_{2n}(t), &t \ge s/2 \end{cases} \]
and
\[ H_2(s,k)(t) = \begin{cases} \nu_{2n}((2t)/(2-s)), &t \le 1/2,\\
\fF_\gamma(2t-1) \nu_{2n}((2t)/(2-s)), &1/2 \le t \le 1 - s/2,\\
\gamma(2t-1), &t \ge 1 - s/2. \end{cases} \]
Straightforward estimates of the expressions above complete the proof.
\qed

This completes the proof of Theorem \ref{theo:homotosur}.
For later use, we want a geometric understanding
of what this tells us about $H^2(\cL_{\pm 1})$.

Recall that $H_2(\cI_{\pm 1}; \ZZ) = \pi_2(\cI_{\pm 1}) =
\pi_2(\Omega\Ss^3) = \ZZ$.
Since an element of $\Omega\Ss^3$ is a function from $\Ss^1$ to $\Ss^3$,
a map $\alpha: \Ss^2 \to \Omega\Ss^3$ can be reinterpreted
as a map $\hat\alpha: \Ss^2 \times \Ss^1 \to \Ss^3$.
The identification between $\pi_2(\Omega\Ss^3)$ and $\ZZ$
takes such a map $\alpha$ to the degree of $\hat\alpha$.

Similarly, let $M$ be a closed oriented surface
and consider a map $\beta: M \to \cI_{+1}$.
Let $\hat\beta: M \times \Ss^1 \to \Ss^3$
be defined by $\hat\beta(p,t) = \tilde\fF_{\beta(p)}(t)$.
Define $\xx: H_2(\cI_{+1}; \ZZ) \to \ZZ$ by
$\xx(\beta) = \deg(\hat\beta)$:
the map $\xx$ provides the identification $H_2(\cI_{+1};\ZZ) = \ZZ$
and is a generator of $H^2(\cI_{+1}; \ZZ)$.
A similar construction defines $\xx \in H^2(\cI_{-1}; \ZZ)$.
The inclusion $\cL_{\pm 1} \subset \cI_{\pm 1}$
defines $\xx \in H^2(\cL_{\pm 1}; \ZZ)$.

As we shall see later, a function $f: K \to \cL_I \subset \cI_I$
may be homotopic to a constant in $\cI_I$
but not in $\cL_I$.
The following proposition shows that this changes if we add loops.

\begin{prop}
\label{prop:easyloop}
Let $K$ be a compact set and 
let $f: K \to \cL_I \subset \cI_I$ a continuous function.
If $f$ is homotopic to a constant in $\cI_I$
then, for any $n > 0$,
$\nu_n \ast f$ is homotopic to a constant in $\cL_I$.
\end{prop}

{\nobf Proof: }
Let $H: K \times [0,1] \to \cI_I$ be a homotopy
with $H(\cdot,0) = f$, $H(\cdot,1)$ constant.
By Lemma \ref{lemma:YtoX}, for sufficiently large $n$,
say $n > N$, the image of $F_n \circ H$ is contained in $\cL_I$.
This implies that $F_n \circ f$ ($n > N$)
is homotopic in $\cL_I$ to a constant.
From Lemma \ref{lemma:Fandp}, $\nu_n \ast f$
is homotopic to $F_n \circ f$ in $\cL_I$ and therefore
the proposition is proved for large $n$.

On the other hand, as we observed in Figure \ref{fig:addloopl},
one loop can be converted to three loops.
The interval $[0,1/2]$ counts as a loop in $\nu_1 \ast f$
and therefore $\nu_1 \ast f$ is homotopic in $\cL_I$ to $\nu_3 \ast f$.
More generally, $\nu_n \ast f$ is homotopic to $\nu_{n+2} \ast f$,
completing the proof.
\qed

\section{Flowers}

\label{sect:flower}

Recall that a curve $\gamma \in \cL_{(-1)^{(k+1)}}$ is a {\sl flower}
of $2k+1$ petals if
there exist $0 = t_0 < t_1 < t_2 < \cdots < t_{2k} < t_{2k+1} = 1$ and
$0 = \theta_0 < \theta_1 < \theta_2 < \cdots < \theta_{2k} < \theta_{2k+1}
= \theta_M$ such that:
\begin{enumerate}
\item{$\gamma(t_1) = \gamma(t_2) = \cdots = \gamma(t_{2k}) = e_1$;}
\item{the only self-intersections of the curve $\gamma$
are of the form $\gamma(t_i) = \gamma(t_j)$;}
\item{the argument of $(x_{i,2}, x_{i,3})$ is $\theta_i$,
where $(0, x_{i,2}, x_{i,3}) = (-1)^i \gamma'(t_i)$.}
\end{enumerate}

For $i = 0, \ldots, 2k$,
it follows from Section \ref{sect:convex} that the restrictions
$\gamma|_{[t_i,t_{i+1}]}$ (the petals) are convex curves.
Let $\cF_{2k} \subset \cL_{(-1)^{(k+1)}}$ be the set of flowers
with $2k+1$ petals.

\begin{lemma}
\label{lemma:flower}
The subset $\cF_{2k} \subset \cL_{(-1)^{(k+1)}}$ is closed.
Also, there is an open neighborhood $\cA_{2k}$ of $\cF_{2k}$
and smooth function $\psi_{2k}: \cA_{2k} \to \RR^{2k}$
such that $0$ is a regular value and $\cF_{2k} = \psi_{2k}^{-1}(\{0\})$.
\end{lemma}

{\nobf Proof: }
We first prove that the sets $\cF_{2k}$ are closed.
Since the region near $e_1$ is taken care of by definition,
all we have to check is that no self-tangencies
within one petal or beween different petals
exist in limit cases of flowers.
Within each petal, a self tangency contradicts the convexity of the petal.
The convexity of petals also implies that the image under $\gamma$
of the interval $(t_i, t_{i+1})$ is contained in the open region defined by
$\langle v, \bfn(t_i) \rangle > 0$, $\langle v, \bfn(t_{i+1}) \rangle > 0$.
Notice that these regions are disjoint and removed from each other
except at the points $\pm e_1$.
Thus, one petal can not touch another petal and therefore $\cF_{2k}$ is closed.

\begin{figure}[ht]
\begin{center}
\psfrag{e1}{$e_1$}
\psfrag{gt1}{$\gamma(\hat t_1)$}
\psfrag{gt2}{$\gamma(\hat t_2)$}
\epsfig{height=30mm,file=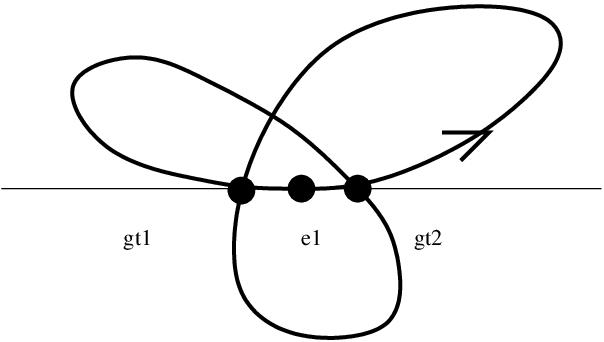}
\end{center}
\caption{A curve $\gamma$ near a flower with $3$ petals.}
\label{fig:almosttrefoil}
\end{figure}

Near a flower $\gamma_0$, curves $\gamma$ will intersect
the large circle through $e_1$ and $e_2$ transversally
(see Figure \ref{fig:almosttrefoil}).
Let $\eta_0(t) = \langle \gamma(t), e_3 \rangle$:
we have $2k$ solutions $\hat t_i \approx t_i$ to $\eta_0(\hat t) = 0$
(where $\gamma_0(t_i) = e_1$).
Define
\[ \psi_{2k}(\gamma) = \left(
\langle \gamma(\hat t_1), e_2 \rangle,
\langle \gamma(\hat t_2), e_2 \rangle, \ldots,
\langle \gamma(\hat t_{2k}), e_2 \rangle \right). \]
Clearly, $\psi_{2k}(\gamma) = 0$ if and only if $\gamma$ is a flower.
The regularity of the value $0$ follows from the fact that the curves
are transversal to the horizontal plane at $\hat t_i$,
completing the proof.
\qed

The sets $\cF_{2k}$ are contractible:
this follows from the fact that the space of possible petals
(convex curves) is contractible (\cite{ShapiroM}) for every choice of 
$t_1, \ldots, t_{2k}$ and $\theta_1, \ldots, \theta_{2k}$.
We shall not use this fact in this paper; it will be proved in \cite{S2}.

Counting intersections with $\cF_{2k}$ defines an element
$\ff_{2k}$ in $H^{2k}(\cL_{(-1)^{(k+1)}}; \ZZ)$.
Since the sets $\cF_{2k}$ are disjoint and the normal bundle is trivial
we have $\ff_{2k}\ff_{2k'} = 0$ both for $k \ne k'$ and $k = k'$.
Also, the degree of $\hat\beta: M \times \Ss^1 \to \Ss^3$
(as in the definition of $\xx$)
can be computed at an element $z \in \Ss^3$ with
\[ \Pi(z) = \begin{pmatrix} -1 & 0 & 0 \\
0 & 1 & 0 \\ 0 & 0 & -1 \end{pmatrix}; \]
since no flower ever passes through the point $-e_1$,
we have $\xx\ff_{2k} = 0$.
We still have to prove that $\ff_{2k} \ne 0$.

\section{Construction of $\bfg_2: \Ss^2 \to \cL_{+1}$}

\label{sect:k=1}

It is probably good to begin by recalling Little's proof that
$\nu_2$ (a circle drawn twice)
and $\nu_4$ (a circle drawn four times)
are in the same connected component of $\cL_I$.
Figure \ref{fig:2to4} below illustrates this.
Initially perturb your curve in order to have three self-intersection
points forming approximately an equilateral triangle.
Pull out the ``petals'' to obtain a flower with three petals,
as in the third figure.
Pull the petals even further so that you have a curvilinear triangle
with loops at the three vertices.
The passage from the fourth to the fifth figure is the only one
where it is important to recall that we are in a sphere,
not in the plane:
one way to think of this is that the triangle became large
and the bulk of the sphere passed through the triangle.
Now it is merely a matter of bringing the three loops together
and making the curve round again.

\begin{figure}[ht]
\psfrag{x2}{$\nu_2$}
\psfrag{x4}{$\nu_4$}
\begin{center}
\epsfig{height=18mm,file=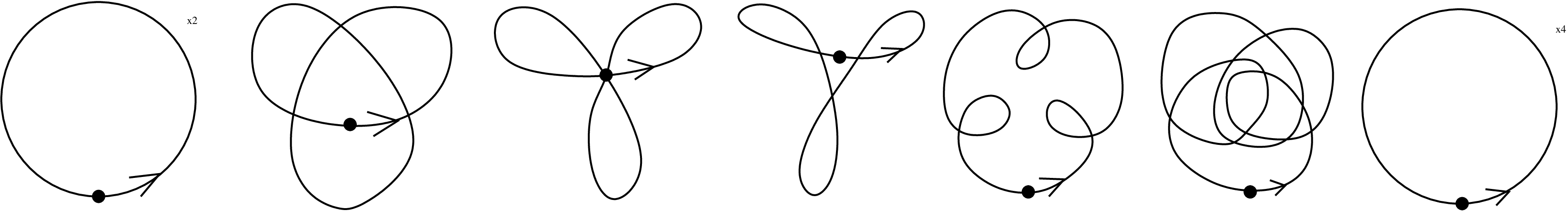}
\end{center}
\caption{A path from $\nu_2$ to $\nu_4$.}
\label{fig:2to4}
\end{figure}

\begin{figure}[ht]
\begin{center}
\epsfig{height=25mm,file=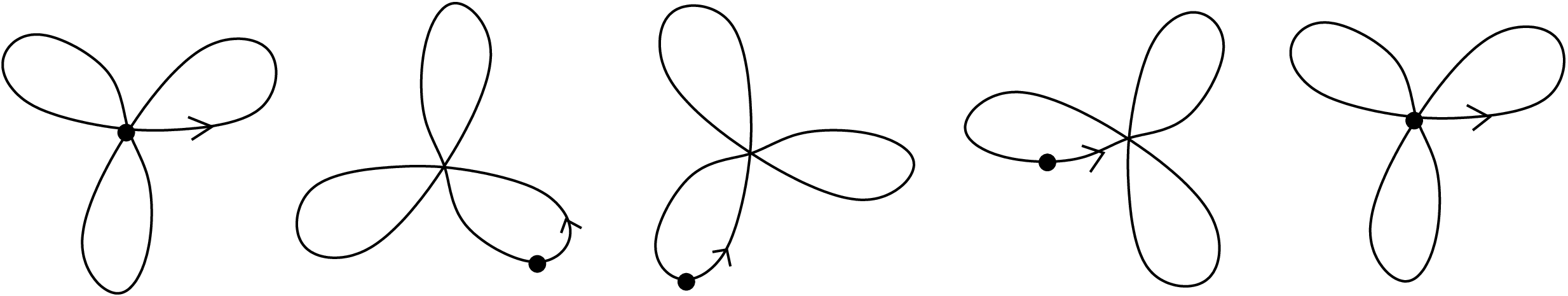}
\end{center}
\caption{The image of a circle under $\bfg_{+,2}$.}
\label{fig:otherline}
\end{figure}

In the construction of this path
there was one important arbitrary choice:
the position of the base point, or, equivalently,
the orientation of the triangle.
From either point of view, the construction can be turned
producing a continuous family indexed by $\Ss^1$ of such paths.
Since the endpoints of all paths coincide,
this is equivalent to constructing a map $\bfg_{+,2}: \Ss^2 \to \cL_{+1}$,
with the two poles taken to $\nu_2$ and $\nu_4$,
meridians (from one pole to the other) corresponding to paths like the one
in Figure \ref{fig:2to4}
and parallel circles being taken to paths obtained by rotating
the curve an in Figure \ref{fig:otherline}.

\begin{figure}[ht]
\begin{center}
\epsfig{height=80mm,file=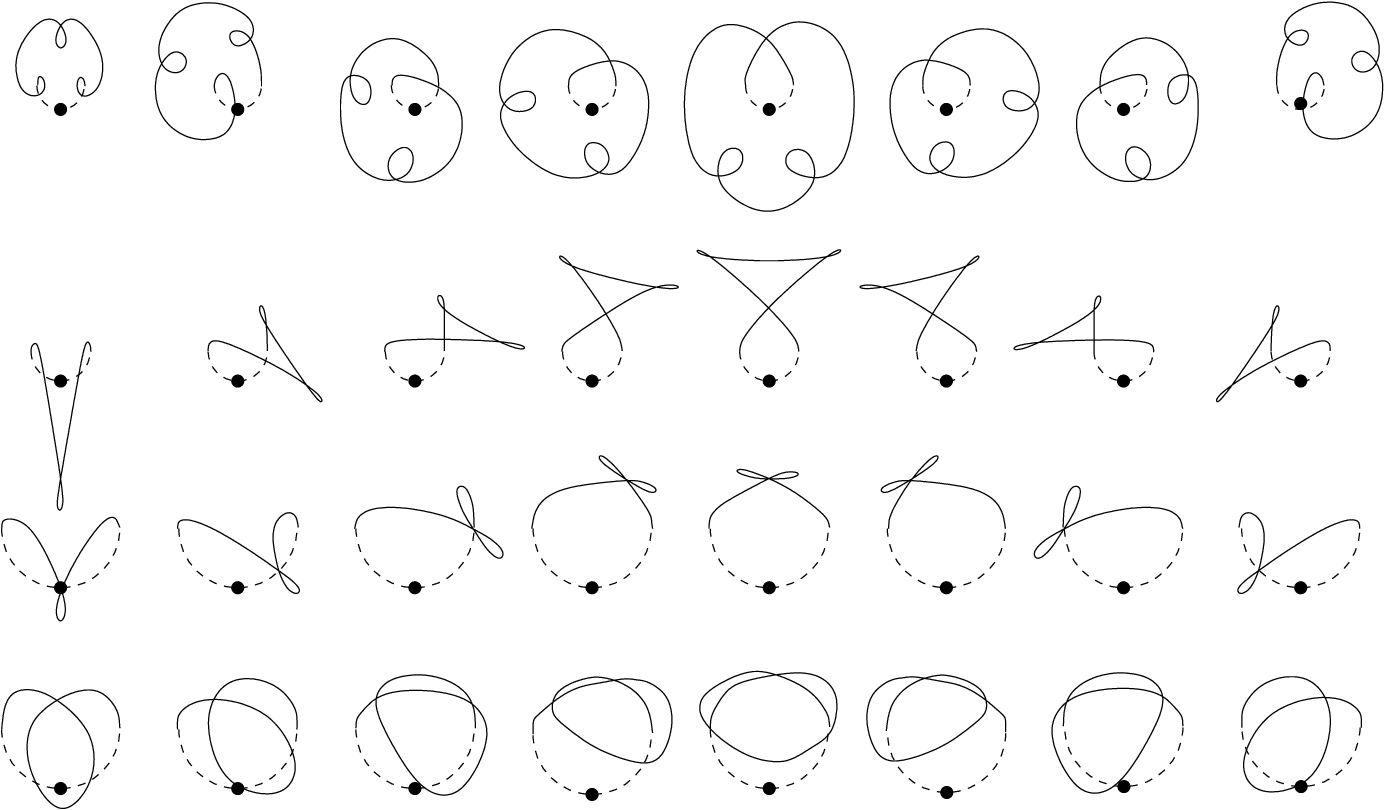}
\end{center}
\caption{The function $\bfg_{+,2}: \Ss^2 \to \cL_{+1}$.}
\label{fig:bfgp2}
\end{figure}

The whole construction is illustrated in Figure \ref{fig:bfgp2}.
The leftmost and rightmost columns are adjacent,
the north pole $\nu_4$ is at the top and the south pole $\nu_2$
is at the bottom (as in most world maps).
In the transition between the first and second rows
most of the curve passed around the back of the sphere.

As an alternative to this figure,
we provide a formula for $\bfg_{+,2}$.
Let $\alpha: [0,1] \times [0,1] \to SO(3)$
be defined by
\[ \alpha(s,t) = \begin{pmatrix}
\sin \pi s \cos 2\pi t & -\sin 2\pi t & -\cos \pi s \cos 2\pi t \\
\sin \pi s \sin 2\pi t &  \cos 2\pi t & -\cos \pi s \sin 2\pi t \\
\cos \pi s & 0 & \sin \pi s \end{pmatrix} \]
and define
\[ \gamma_s(t) = \frac{\sqrt{2}}{2}\;\alpha(s,t)
\begin{pmatrix} 1 \\ \cos 6\pi t \\ \sin 6\pi t \end{pmatrix}. \]
The curve $\gamma_0$ is a circle drawn $4$ times and
the curve $\gamma_1$ is a circle drawn $2$ times.
A tedious but straightforward computation verifies that
the curves $\gamma_s$ are closed and locally convex
and therefore belong to $\cL_{per}$.

Let $\Gamma(s,t) = \fF_{\gamma_s}(t)$:
it easy to verify that 
\[ \Gamma(s,t+(1/3)) =
\begin{pmatrix} -1/2 & -\sqrt{3}/2 & 0 \\
\sqrt{3}/2 & -1/2 & 0 \\ 0 & 0 & 1 \end{pmatrix} \Gamma(s,t) \]
for all $s$ and $t$.
Finally, let $\bfg_{+,2}: [0,1] \times [0,1] \to \cL_{+1}$ be defined by
\[ \bfg_{+,2}(s_1, s_2)(t) =
(\Gamma(s_2,s_1/3))^{-1} \Gamma(s_2,t + (s_1/3)) e_1. \]
If $s_2 = 0$ or $1$, the value of $s_1$ is irrelevant 
for the value of $\bfg_{+,2}$.
Also, $\bfg_{+,2}(0,s_2) = \bfg_{+,2}(1,s_2)$ for all $s_2$.
Performing these identifications, the domain of $\bfg_{+,2}$ becomes
the sphere $\Ss^2$, as required.

It follows easily either from Figure \ref{fig:bfgp2}
or from the formulas that $\xx(\bfg_{+,2}) = 1$,
i.e., that the degree of $\hat\bfg_{+,2}: \Ss^2 \times \Ss^1 \to \Ss^3$,
$\hat\bfg_{+,2}(p,t) = \tilde\fF_{\bfg_{+,2}(p)}(t)$, equals $1$.
This can be seen, for instance,
by computing preimages of some $z \in \Ss^3$.
Thus, $\bfg_{+,2}$ is a generator of $\pi_2(\cI_{+1})$.

It is again clear from Figure \ref{fig:bfgp2} that $\bfg_{+,2}$
intersects $\cF_2$ precisely once
(the flower is in the third row, first column),
and the intersection is transversal,
and therefore $|\ff_2(\bfg_{+,2})| = 1$.
We were not too careful about the orientation of $\ff_2$
in Section \ref{sect:flower} so we decree now that $\ff_2(\bfg_{+,2}) = 1$.
This can likewise be checked for the formula by a long
and tedious computation which we skip.
Either way, $\ff_2 \ne 0$.

Consider the function $\nu_2 \ast \bfg_{+,2}$.
We have $\xx(\nu_2 \ast \bfg_{+,2}) = 1$
(adding these loops does not change the degree)
and $\ff_2(\nu_2 \ast \bfg_{+,2}) = 0$
(there are no flowers in the image of $\nu_2 \ast \bfg_{+,2}$
since no flower starts with two loops).
In particular, the maps $\bfg_{+,2}$ and $\nu_2 \ast \bfg_{+,2}$
are not homotopic in $\cL_{+1}$.
On the other hand, from Lemma \ref{lemma:addloopi},
the maps $\bfg_{+,2}$ and $\nu_2 \ast \bfg_{+,2}$
are homotopic in $\cI_{+1}$.

We will now consider the group $\pi_2(\cL_{+1})$
but before we do so we must say a few words about base points.
Recall that given two base points $p_1$ and $p_2$,
the two homotopy groups $\pi_2(\cL_{+1};p_1)$ and $\pi_2(\cL_{+1};p_2)$
are identified via a homotopy class of paths from $p_1$ to $p_2$.
We prove in \cite{S2} that $\cL_{+1}$ and $\cL_{-1,n}$
are simply connected and therefore the identification is natural;
the idea here is not, however, to use these results.
We shall therefore select $\nu_2$ as a base point for $\cL_{+1}$
and $\nu_3$ as a base point for $\cL_{-1}$.

\begin{figure}[ht]
\begin{center}
\psfrag{bfgp2}{$\bfg_{+,2}$}
\psfrag{nu2bfgp2}{$\nu_2 \ast \bfg_{+,2}$}
\psfrag{alpha}{$\alpha$}
\epsfig{height=30mm,file=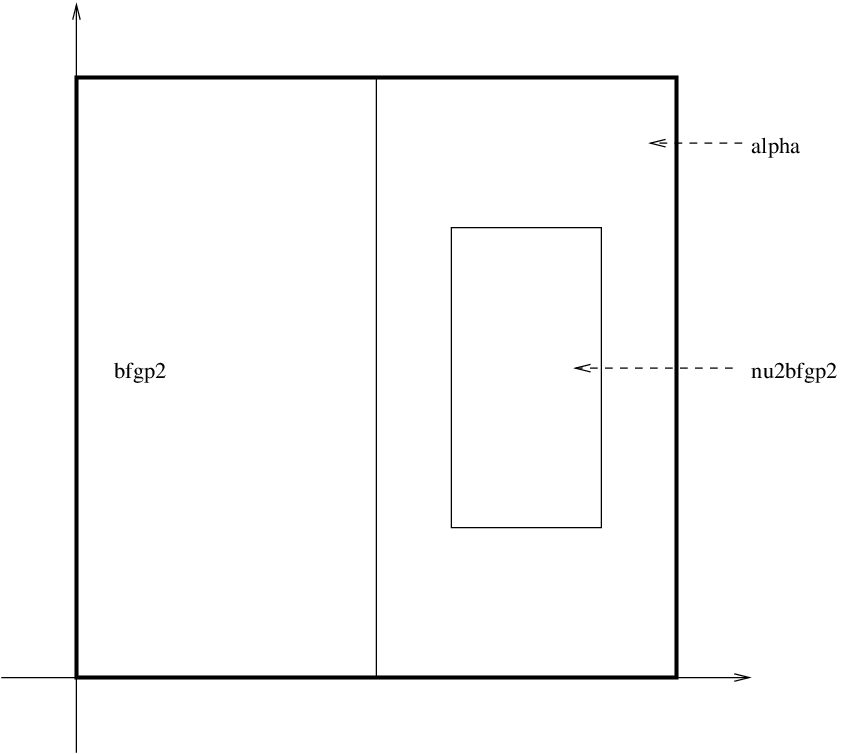}
\end{center}
\caption{The construction of $\bfg_{2}$.}
\label{fig:bfg2}
\end{figure}

Finally, consider the difference
$\bfg_2 = \bfg_{+,2} - (\nu_2 \ast \bfg_{+,2}): \Ss^2 \to \cL_{+1}$.
More precisely, consider $\bfg_{+,2}$ as a function from $[0,1]^2$
to $\cL_{+1}$ with $\bfg_{+,2}(p) = \nu_2$ for $p \in \partial([0,1]^2)$.
Let $\alpha: [0,1] \to \cL_{+1}$ be a path from $\nu_2$ to $\nu_4$ and define
$\bfg_2: [0,1]^2 \to \cL_{+1}$ by
\[ \bfg_2(x,y) = \begin{cases} \bfg_{+,2}(2x,y), & 0 \le x \le \frac{1}{2}; \\
(\nu_2 \ast \bfg_{+,2})(4x - \frac{5}{2}, 2y - \frac{1}{2}), &
\max(2|x-\frac{3}{4}|,|y-\frac{1}{2}|) \le \frac{1}{4}; \\
\alpha(1-4 \max(2|x-\frac{3}{4}|,|y-\frac{1}{2}|) ), &
\frac{1}{4} \le \max(2|x-\frac{3}{4}|,|y-\frac{1}{2}|) 
 \le \frac{1}{2}. \end{cases} \]
This construction is sketched in Figure \ref{fig:bfg2}.
From Lemma \ref{lemma:addloopi},
this map is homotopic to a constant in $\cI_{+1}$.
The path $\alpha$ can be chosen so as to avoid the set $\cF_2$
and therefore the image of $\bfg_2$ intersects $\cF_2$
transversally and exactly once and we have $\ff_2(\bfg_2) = 1$.
This implies that 
neither $\bfg_2$ nor any positive multiple of it
is homotopic to a constant in $\cL_{+1}$.
We therefore have a copy of $\ZZ^2$ contained in $\pi_2(\cL_{+1})$.
This completes the proof of the following result,
closely related to the case $k=1$ of Theorem \ref{theo:Gcohomology}.

\begin{theo}
\label{theo:k=1}
The maps $\bfg_{+,2}$ and $\bfg_2$ span a copy of $\ZZ^2$
contained in $\pi_2(\cL_{+1})$.
The elements $\xx$ and $\ff_2$ span a copy of $\ZZ^2$
contained in $H^2(\cL_{+1}; \ZZ)$.
\end{theo}

\section{The homotopy lifting property}

\label{sect:homolift}

Recall that $\phi: \cL \to \Ss^3$ takes $\gamma$
to $\tilde\fF_{\gamma}(1)$.
Let $\DD^2$ be the closed unit disk,
$X_1 = \DD^2 \times [0,1]$ and
$X_2 = \DD^2 \times \{0\} \cup \Ss^1 \times [0,1]$.
We construct functions $\bfh: X_1 \to \Ss^3$
and $\tilde\bfh_0: X_2 \to \cL$ with
$\phi \circ \tilde\bfh_0 = \bfh|_{X_2}$.
We then prove
that there is no continuous function $\tilde\bfh: X_1 \to \cL$
with $\phi \circ \tilde\bfh = \bfh$,
thus proving that the homotopy lifting property does not hold.


For $c \in (0,+\infty)$,
let $\nu_c \in \cL$ be defined by $\nu_c(t) = \nu_1(ct)$
so that $\phi(\nu_c) = \tilde\fF_{\nu_c}(1) = \tilde\fF_{\nu_1}(c)$,
\[ \fF_{\nu_1}(c) =
\frac{1}{2} \; \begin{pmatrix}
1 + \cos(2\pi c) & - \sqrt{2} \sin(2\pi c) & 1 - \cos(2\pi c) \\
\sqrt{2} \sin(2\pi c) & 2 \cos(2\pi c) & -\sqrt{2} \sin(2\pi c) \\
1 - \cos(2\pi c) & \sqrt{2} \sin(2\pi c) & 1 + \cos(2\pi c)
\end{pmatrix}. \]
Define $\bfh: X_1 = \DD^2 \times [0,1] \to \Ss^3$ by
$\bfh(p,s) = \phi(\nu_{4-2s}) = \tilde\fF_{\nu_4}(1-s/2)$;
notice that $\bfh(p,0) = \bfh(p,1) = 1$ for all $p \in \DD^2$.

\begin{figure}[ht]
\begin{center}
\psfrag{s=0}{$s=0$}
\psfrag{s=1}{$s=1$}
\psfrag{nu2}{$\nu_2$}
\psfrag{nu4}{$\nu_4$}
\psfrag{nus}{$\nu_{4-2s}$}
\psfrag{flower}{flower}
\psfrag{bfg2}{$\bfg_{+,2}(\pi_{\DD^2,\Ss^2}(p))$}
\epsfig{height=50mm,file=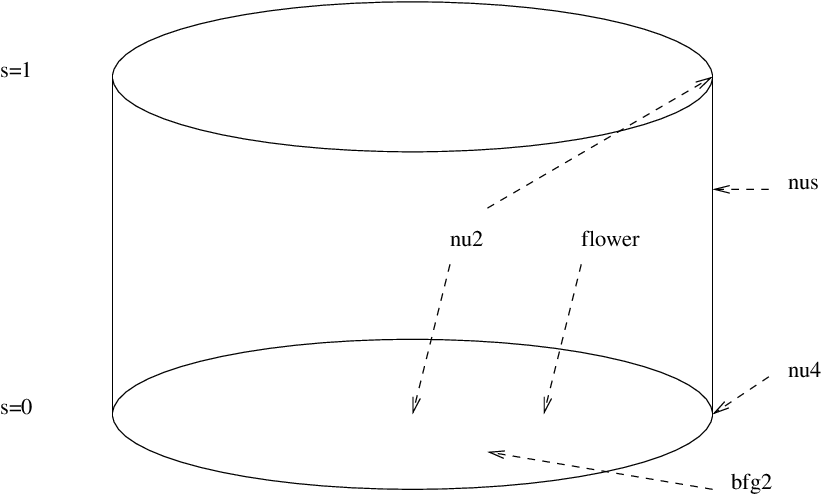}
\end{center}
\caption{The function $\tilde\bfh_0:
X_2 = \DD^2 \times \{0\} \cup \Ss^1 \times [0,1] \to \cL$.}
\label{fig:tbfh0}
\end{figure}

Let $\pi_{\DD^2,\Ss^2}: \DD^2 \to \Ss^2$ be the function
that wraps the sphere by taking
the boundary of $\DD$ to the north pole of $\Ss^2$,
other points of $\Ss^2$ having one transversal preimage.
Define $\tilde\bfh_0(p,0) = (\bfg_{+,2} \circ \pi_{\DD^2,\Ss^2})(p)$
for all $p \in \DD^2$; notice that 
$\tilde\bfh_0(p,0) = \nu_4$ for $p \in \Ss^1 = \partial\DD^2$.
Finally, for $p \in \Ss^1$,
define $\tilde\bfh_0(p,s) = \nu_{4-2s}$ (see Figure \ref{fig:tbfh0}).

Assume by contradiction that $\tilde\bfh: X_1 \to \cL$
satisfies $\phi \circ \tilde\bfh = \bfh$.
We construct $\bfh_3: X_1 \to \cL_{+1}$ by completing
the locally convex curves $\tilde\bfh(p,s)$ with
a parametrized arc of $\cC_0$. More precisely,
\[ (\bfh_3(p,s))(t) = \begin{cases}
(\tilde\bfh(p,s))\left(\frac{t}{1-s/2}\right), & 0 \le t \le 1 - s/2, \\
\nu_4(t), & 1 - s/2 \le t \le 1; \end{cases} \]
in particular, $\bfh_3(p,s) = \nu_4$ for $p \in \Ss^1$.
The construction above guarantees the continuity of $\bfh_3$.

We now consider the 2-cycle $\bfh_2 = \bfh_3|_{\partial X_1}$
and its product with $\ff_2$.
In other words, we count flowers in the image of the boundary.
There is a unique flower at the image of the bottom $\DD^2 \times \{0\}$:
this intersection with $\cF_2$ is transversal.
Since $\bfh_2(p,s) = \nu_4$ for all $(p,s) \in \Ss^1 \times [0,1]$
there are no flowers on the sides.
Finally, $\bfh_2(p,1)(t) = \nu_4(t)$ for all $t \ge 1/2$:
curves on the top finish with two turns around $\cC_0$
and are therefore definitely not flowers.
This means that the product of $\bfh_2$ with $\ff_2$ is not zero
and therefore $\bfh_2 \ne 0 \in \pi_2(\cL_{+1})$,
contradicting the existence of $\bfh_3$.

\section{Construction of $\bfg_{2k}: {\Ss^2}^k \to \cL_{(-1)^{(k+1)}}$ and \\
proof of Theorem \ref{theo:Gcohomology}}

\label{sect:Gcohomology}

There is a dashed arc around the base point in each curve
in Figure \ref{fig:bfgp2}.
The dashed arc remains unchanged during the entire process.
We show that a minor modification of $\bfg_{+,2}$ can be constructed
so that this dashed arc is a circle minus a small gap,
i.e., changes are restricted to a small interval.

Take any function $f: K \to \cL_I$, $K$ compact.
For sufficiently small $\epsilon_1 > 0$, the arcs
$\gamma|_{[0,\epsilon_1]}$ and $\gamma|_{[1-\epsilon_1,1]}$
are convex for any $\gamma = f(p)$, $p \in K$.
More, for sufficiently small $\epsilon_2 > 0$ the arcs
$\nu_1|_{[0,\epsilon_2]}$ and $\nu_1|_{[1-\epsilon_2,1]}$
can be inserted in $\gamma$ without damaging convexity.
Thus a homotopy $H: [0,1] \times K \to \cL_I$, $H(0,p) = f(p)$,
changes the curves only in a small neighborhood of the base point
and at the end of the homotopy we have
$\gamma(t) = \nu_1(t)$ for all $\gamma = H(1,p)$, $p \in K$
and for all $t \in [0,\epsilon_2] \cup [1-\epsilon_2,1]$.

Let $R$ be an upper triangular matrix with positive diagonal
and $\gamma \in \cL_I$ a locally convex curve.
The curve $\alpha: [0,1] \to \RR^3$, $\alpha(t) = R^{-1} \gamma(t)$, satisfies
$\det(\alpha(t), \alpha'(t), \alpha''(t)) > 0$ for all $t$.
The curve $\gamma^R: [0,1] \to \Ss^2$,
$\gamma^R(t) = \alpha(t)/|\alpha(t)|$ also satisfies
$\det(\gamma^R(t), (\gamma^R)'(t), (\gamma^R)''(t)) > 0$ for all $t$
and therefore is locally convex.
The group of matrices of the form
\[ R = \begin{pmatrix} a^{-1} & a^{-1} b & a^{-1} b^2/2 \\
0 & 1 & b \\ 0 & 0 & a \end{pmatrix}, \quad a > 0, \]
takes the cone $y^2 = 2xz$ onto itself.
Thus, for $R$ as above, the small arc of the circle $\cC_0$ around $e_1$
common to all curves $\gamma$, $\gamma = f(p)$, $p \in K$,
is taken to another arc of $\cC_0$ common to all curves $\gamma^R$,
with another parametrization different from $\nu_1$
but common to all curves.
If $a$ is taken to be a large positive number,
the arc will be arbitrarily large, $\cC_0$ minus a small gap;
an appropriate choice of $b$ allows us to position
that gap anywhere along $\cC_0$ away from $e_1$.
A reparametrization allows us to assume that
$\gamma(t) = \nu_1(t)$ except in a small interval $I \subset [0,1]$.
Notice that this construction preserves the fact that there is
only one intersection with $\cF_2$, and this intersection is transversal.

We are now ready to construct $\bfg_{2k}$ recursively
from $\bfg_2$ and $\bfg_{2k-2}$.
Assume by induction that $\bfg_{2k-2}$ intersects the manifold $\cF_{2k-2}$
transversally and exactly once so that $\ff_{2k-2}(\bfg_{2k-2}) = 1$
and that $\nu_n \ast \bfg_{2k-2}$ is homotopic to a constant for any $n > 0$.
We first construct $\bfg_{+,2k}$ with domain
$\Ss^{2k} = \DD^2 \times \Ss^{2k-2} \cup \Ss^1 \times \DD^{2k-1}$.

\begin{figure}[ht]
\begin{center}
\psfrag{g1}{$g_1$}
\psfrag{g2}{$g_2$}
\epsfig{height=50mm,file=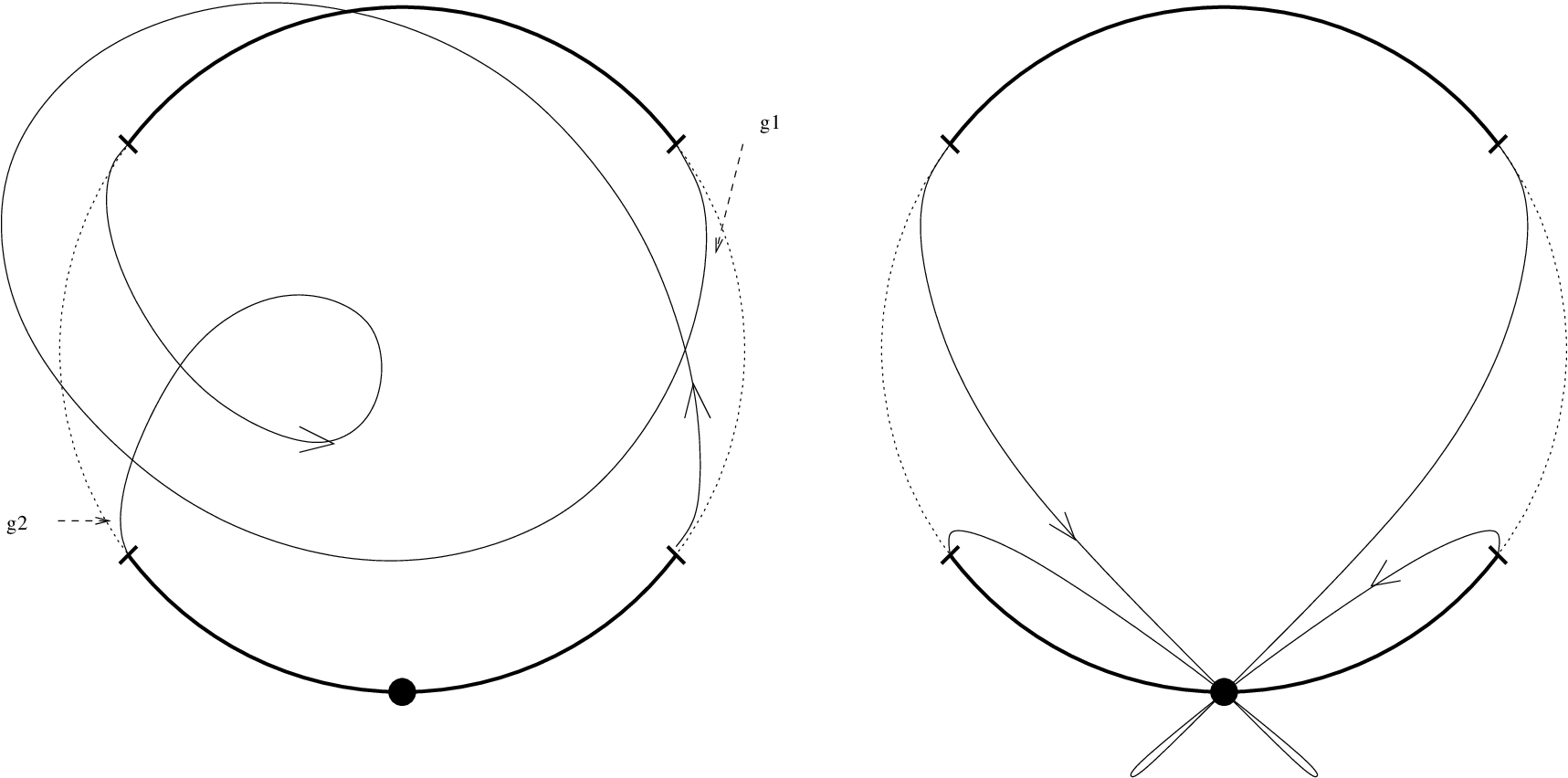}
\end{center}
\caption{Two curves in the image of $\bfg_{+,2k}$; the second is a flower.}
\label{fig:bfg2k}
\end{figure}

Let $I_1 = [1/6,2/6]$, $I_2 = [4/6,5/6]$.
Define functions
$g_1: \DD^2 \to \cL_{+1}$ and $g_2: \Ss^{2k-2} \to \cL_{(-1)^k}$ with:
\begin{enumerate}[(a)]
\item{$g_1(p)(t) = \nu_1(t)$ for all $p \in \DD^2$,
$t \in [0,1] \smallsetminus I_1$;}
\item{$g_1(p)$ is a reparametrization of $\nu_4$
for $p \in \Ss^1 = \partial\DD^2$;}
\item{$g_1$ intersects $\cF_2$ transversally and exactly once
at $p_1^\cF \in \int(\DD^2)$;}
\item{$g_2(p)(t) = \nu_1(t)$ for all $p \in \Ss^{2k-2}$,
$t \in [0,1] \smallsetminus I_2$;}
\item{$g_2$ intersects $\cF_{2k-2}$ transversally and exactly once
at $p_2^\cF \in \Ss^{2k-2}$;}
\item{$\nu_n \ast g_2$ is homotopic to a constant in $\cL_I$
for $n > 0$.}
\end{enumerate}
The function $g_1$ is obtained from $\bfg_{+,2} \circ \pi_{\DD^2,\Ss^2}$
via the above construction.
Similarly, the function $g_2$ is obtained from $\bfg_{2k-2}$
by the same construction. For $(p_1,p_2) \in \DD^2 \times \Ss^{2k-2}$
define
\[ \bfg_{+,2k}(p_1,p_2)(t) = \begin{cases}
g_1(p_1)(t),&t \in I_1, \\
g_2(p_2)(t),&t \in I_2, \\
\nu_1(t),& \textrm{otherwise.} \end{cases} \]
As in Figure \ref{fig:bfg2k}, we can say that
each $g_j$ is responsible for the interval $I_j$.
Notice that for $p_1 \in \Ss^1$, $\bfg_{+,2k}(p_1,p_2)$
is a reparametrization of $\nu_1 \ast (\nu_2 \ast g_2(p_2))$
(the reparametrization is independent of $p_1$).
Let $g_3: \DD^{2k-1} \to \cL_{(-1)^k}$
be a continuous map with
$g_3(p_2) = \nu_2 \ast g_2(p_2)$ for all $p_2 \in \Ss^{2k-2}$.
Up to the above mentioned reparametrization,
for $(p_1,p_2) \in \Ss^1 \times \DD^{2k-1}$ define
$\bfg_{+,2k}(p_1,p_2) = \nu_1 \ast g_3(p_2)$.

As in Figure \ref{fig:bfg2k},
the only intersection of $\bfg_{+,2k}$ with $\cF_{2k}$
in $\DD^2 \times \Ss^{2k-2}$ is in $(p_1^\cF, p_2^\cF)$;
this intersection is transversal.
There are no intersections of $\bfg_{+,2k}$ with $\cF_{2k}$
in $\Ss^1 \times \DD^{2k-1}$ for curves there are of the form
$\nu_1 \ast \textrm{(something)}$.

Finally, construct $\bfg_{2k} = \bfg_{+,2k} - (\nu_2 \ast \bfg_{+,2k})$,
as for $\bfg_2$ in Section \ref{sect:k=1}.
By Lemma \ref{lemma:addloopi} and Proposition \ref{prop:easyloop}, 
$\bfg_{2k}$ is homotopic to a constant in $\cI_I$ or, equivalently,
$\nu_2 \ast \bfg_{2k}$ is homotopic to a constant in $\cL_I$.
There are no flowers in the image of $\nu_2 \ast \bfg_{+,2k}$
and therefore $\ff_{2k}(\bfg_{2k}) = 1$.
This completes the inductive construction of $\bfg_{2k}$,
proves that $\ff_{2k} \ne 0$ and completes the proof
of Theorem \ref{theo:Gcohomology}.
We sum up some of our other conclusions as another theorem.

\begin{theo}
\label{theo:pik}
Let $k \ge 1$.
Consider the inclusion $i: \cL_{(-1)^{(k+1)}} \to \cI_{(-1)^{(k+1)}}$
and the induced map $\pi_{2k}(i): 
\pi_{2k}(\cL_{(-1)^{(k+1)}}) \to \pi_{2k}(\cI_{(-1)^{(k+1)}})$.
Then the map $\bfg_{2k}$ constructed above spans a copy of $\ZZ$
in $\ker(\pi_{2k}(i))$.
\end{theo}

\section{Final remarks}

\label{sect:conclusion}

In the second paper of this series (\cite{S2}) we
show that connected components of $\cL_I$ are simply connected.
We also show that the inclusion $\cL_{-1,n} \subset \cI_{-1}$
induces an isomorphism between $\pi_2(\cL_{-1,n})$ and $\pi_2(\cI_{-1}) = \ZZ$
and that $\bfg_2$ and $\tilde\bfg_2$ (as in Section \ref{sect:k=1})
actually generate $\pi_2(\cL_{+1}) = \ZZ^2$.
This implies that $H^2(\cL_{+1}; \ZZ) = \ZZ^2$ and
$H^2(\cL_{-1}; \ZZ) = \ZZ$.
In \cite{Saldanha3} we hope to prove that the classes $\xx^n$ and $\ff_{2n}$
are generators of $H^\ast(\cL_{\pm 1})$
and that $\cL_{+1}$ and $\cL_{-1,n}$ have the homotopy type of
$\Omega\Ss^3 \vee \Ss^2 \vee \Ss^6 \vee \Ss^{10} \vee \cdots$ and
$\Omega\Ss^3 \vee \Ss^4 \vee \Ss^8 \vee \Ss^{12} \vee \cdots$,
respectively.

Little's Theorem that convex curves form a separate connected component
can be rephrased as saying that $\cF_0$, the set of flowers with $1$ petal,
obtains a new element in $H^0(\cL_I)$.
From this point of view that result is the case $k = 0$ of
Theorem \ref{theo:Gcohomology}.

The sets $\cL_{+1}$ and $\cL_{-1}$ can naturally be considered
as two instances of a big family of spaces
$\cL_z = \phi^{-1}(\{z\})$, $z \in \Ss^3$
where $\phi: \cL \to \Ss^3$ takes $\gamma$ to $\tilde\fF_\gamma(1)$.
As we saw in Section \ref{sect:homolift},
this map does not satisfy the homotopy lifting property.
These results also imply that Gromov's $h$-principle (\cite{EM}, \cite{Gromov})
fails for $\phi: \cL \to \Ss^3$.
In \cite{SaSha}, on the other hand,
we show that every space $\cL_z$ is homeomorphic
to either $\cL_{+1}$, $\cL_{-1}$ or $\cI_{+1}$.
Since $\cI_{+1}$ is well understood, this leaves out only the two spaces
studied in this paper.

Finally, similar questions can be asked about curves in $\Ss^n$, $n > 2$
($\gamma$ is locally convex if $\det(\gamma(t), \ldots, \gamma^{(n)}(t)) > 0$);
in \cite{SaSha} we show a few results about these spaces.

\bigbreak


\bigskip\bigskip\bigbreak

{

\parindent=0pt
\parskip=0pt
\obeylines

Nicolau C. Saldanha, PUC-Rio
saldanha@puc-rio.br; http://www.mat.puc-rio.br/$\sim$nicolau/



\smallskip

Departamento de Matem\'atica, PUC-Rio
R. Marqu\^es de S. Vicente 225, Rio de Janeiro, RJ 22453-900, Brazil

}

\end{document}